\documentclass[12pt,a4paper,reqno,intlimits]{amsart}
\usepackage{amsfonts, amsmath, amsxtra, amsthm, amssymb, latexsym, mathrsfs, srcltx, dsfont, euscript}

\usepackage{hyperref}
\hypersetup
{  pdftitle={Shown in AR File Information},
   pdfstartview=FitH,   
   bookmarks=true,      
}

\theoremstyle{plain}

\newtheorem*{theorem*}{Theorem}
\newtheorem*{corollary*}{Corollary}
\newtheorem*{lemma*}{Lemma}
\newtheorem*{property*}{Property}
\newtheorem{proposition*}{Proposition}
\newtheorem*{statement*}{Statement}

\newtheorem{theorem}{Theorem}

\newtheorem{statement}{Statement}

\theoremstyle{definition}

\newtheorem*{definition*}{Definition}
\newtheorem*{remark*}{Remark}
\newtheorem*{example*}{Example} 

\newtheorem{definition}{Definition}

\allowdisplaybreaks

\DeclareMathOperator{\Lin}{Lin}%
\DeclareMathOperator{\Clos}{Clos}%

\makeatletter
\renewcommand{\@biblabel}[1]{#1.\hfill}
\makeatother

\renewcommand{\title}[1]{
\begin{center}
{\bfseries #1}
\end{center}}

\renewcommand{\author}[1]
{\begin{center} {\bfseries{\copyright~~~#1}}
\end{center}}

\renewcommand{\address}[1]
{\begin{center}
\tiny
{\scshape{#1}}
\end{center}\medskip}

\begin{document}

\pagenumbering{arabic}

\thispagestyle{empty}

\markboth{{\footnotesize{\it \textbf{Karpenko~I.\,I., Tyshkevich~D.\,L.}}}} {{\footnotesize {\it \textbf{On The Model Of A Skew--Selfadjoint Operator With A Simple Spectrum
On A\ldots}}}}

\medskip
\normalsize

\title{ON THE MODEL OF A SKEW--SELFADJOINT OPERATOR WITH A SIMPLE SPECTRUM
ON A HILBERT QUATERNION MODULE}

\author{Karpenko~I.\,I., Tyshkevich~D.\,L.}

\address{Taurida National V.\,I.~Vernadsky University \\
Department of Mathematics and Informatics \\
Vernadsky ave., 4, Simferopol, 95007, Ukraine \\
e-mail: \textit{dtyshk@inbox.ru, i\_karpenko@inbox.ru}}

\subsection*{Abstract}
In this work we construct the model of a skew--selfadjoint
operator with a simple spectrum acting on a Hilbert quaternion
bimodule. This result is based on the Spectral Theorem for a
skew--selfadjoint operator. In the case of a bounded normal
operator this Spectral Theorem was announced in the report
\cite{kt_KarTyshkSpDecBedl}. The more detailed research was carried
out in the paper \cite{kt_KarTyshkSpDecUchZ}.

In the same paper we pointed out that the given reasonings enable
us to prove corresponding results for a unbounded
skew--selfadjoint operators too. Results of this article
essentially develop the results%
\footnote%
{%
  Note that since the issue of \cite{kt_Visv} investigations of spectral problems
  in \textit{infinite--dimensional} Hilbert (bi)modules were not carried out as far as we know.
  As it seems for us our papers \cite{kt_KarTyshkSpDecUchZ,kt_KarTyshkMatSt}
  are the first in this row after \cite{kt_Visv}.
}%
\ of the paper \cite{kt_Visv}.

\subsection*{Keywords}
Quaternion, Hilbert space, quaternion module, quaternion bimodule, spectral theorem, normal operator,
skew--selfadjoint operator, operator with a simple spectrum, model of a linear operator,
operator of left multiplication by an independent variable.


\subsection*{2010 Mathematics Subject Classification (MSC2010)}
47B15, 47B25, 47B37

\section*{Introduction}

\subsection*{On the skew field of quaternions $\mathbb{H}$ and $\mathbb{R}$--containing subfields of $\mathbb{H}$}
For reader's convenience recall some definitions and facts
concerning quaternions.

The real quaternionic skew-field $\mathbb{H}$ is a
four-dimensional associative division algebra of a rank 4 over
$\mathbb{R}$ with the basis $\{1,i,j,k\}$ and multiplication rules
\begin{equation*}
   \begin{array}{lll}
      i^2=-1 & j^2=-1 & k^2=-1 \\
      ij=k   & jk=i   & ki=j   \\
      ji=-k  & kj=-i  & ik=-j
   \end{array}
\end{equation*}

For any $q\in\mathbb{H}$ there exist unique
$q_0,\,q_1,\,q_2,\,q_3\in\mathbb{R}$ such that
$q=q_0+q_1i+q_2j+q_3k$ (\textit{real representation} of $q$). Also
it is useful to deal with \textit{the vector form} of quaternion
$q=q_0+\vec{q}$ where $\vec{q}=q_1i+q_2j+q_3k$ is the
\textit{vector} or \textit{imaginary part} of $q$ (if $q=\vec{q}$
then $q$ is called to be a \textit{vector} or an
\textit{imaginary} quaternion).

For instance the vector form  for the product of quaternions $q$
and $p$ is nothing more than the well-known formula of
multiplication
$%
  qp=
  q_0p_0-(\vec{q},\vec{p})+
  \big(
       [\vec{q},\vec{p}\,]+p_0\vec{q}+q_0\vec{p}\,
  \big)
$%
. Here $(\cdot,\cdot)$ is the usual scalar product and
$[\cdot,\cdot]$ is the vector product on the three--dimensional
space $\mathbb{R}\langle i,j,k\rangle$ of vector quaternions.
\textit{Conjugate} of $q$ is defined by
$\overline{q}=q_0-q_1i-q_2j-q_3k$. The map $q\rightarrow
\overline{q}$ is an involution on $\mathbb{H}$, and
$q\overline{q}=\overline{q}q=\sum_{t=0}^3q_t^2\in\mathbb{R}$. One
can define the \textit{absolute value} $|q|$ of $q$ by
$|q|=(q\overline{q})^{1/2}$. Thus, we can consider $\mathbb{H}$ as
a \textit{normed algebra}. An imaginary quaternion whose absolute
value equals 1 is called an \textit{imaginary unit}. From this point of view
one can consider the set of complex numbers $\mathbb{C}$  as
a real subalgebra of $\mathbb{H}\colon\mathbb{C}=\mathbb{R}\langle
1,i\rangle$. Inclusion $\mathbb{C}$ in $\mathbb{H}$ allows us to
obtain the \textit{complex} (or \textit{symplectic})
representation of a quaternion. Namely, for $q=q_0+q_1i+q_2j+q_3k$
we have $q=(q_0+q_1i)+(q_2+q_3i)j=z_1+z_2j$ where $z_1,\,z_2\in
\mathbb{C}$.

Often it is useful to restrict not oneself by choice of a concrete
field; e.g. $\mathbb{C}$. In fact $\mathbb{C}$ is just a specimen
of a field in $\mathbb{H}$,  which extends $\mathbb{R}$. Everywhere
in this paper we denote such a field by $\mathbb{F}$. $\mathbb{F}$
is a commutative and associative division algebra over
$\mathbb{R}$; and dimension of $\mathbb{F}$ is greater than 1.
Hence by Frobenius theorem dimension of $\mathbb{F}$ equals 2; and
$\mathbb{F}$ is isomorphic to $\mathbb{C}$.

A field $\mathbb{F}$ is uniquely determined by some nonreal
quaternion. Indeed, let $q \in \mathbb{F}\setminus \mathbb{R}$.
Write down this quaternion in the vector form: $q=q_0+\vec{q}$.
Then we have $\vec{q}=q-q_0 \in \mathbb{F}\setminus\{0\}$. Let
$f:= \frac{1}{|\vec{q}|}\,\vec{q}$. The vector system $\{1,f\}$ is
linear independent. Hence it is a basis of $\mathbb{F}$ as a
two--dimensional $\mathbb{R}$--algebra. These arguments show that
any two nonreal elements of $\mathbb{F}$ have proportional vector
parts, and this condition is necessary and sufficient for
quaternions with the corresponding vector parts to commute. Thus,
we can characterize any subfield of $\mathbb{H}$ satisfying the
above conditions as a set of all quaternions commuting with some
fixed nonreal quaternion. In addition this necessary and
sufficient condition implies possibility for a corresponding
imaginary unit $f$ to be determined by $\mathbb{F}$ up to $\pm1$.

Consider now  $\mathbb{H}$ as a real Euclidean space with
dimension 4. Choose any normed quaternion $\phi$ to be orthogonal
to $1,f$. Then we have $\phi^2=(\phi,\phi)=-1$. The quaternion
$f\phi\,(=[f,\phi])$ is orthogonal to $1,f,\phi$; and
$f\phi=[f,\phi]=-[\phi,f]=-\phi f$. Hence $(f\phi)^2=-1$. Thus,
the system $\{1,f,\phi,f\phi\}$ is an $\mathbb{R}$--basis of
$\mathbb{H}$; and this basis consists of 1 and three imaginary
units. This fact allows us to define the unique decomposition
$%
  q=q_0+\widetilde{q}_1f+\widetilde{q}_2\phi+\widetilde{q}_3f\phi
$%
\ which implies an analogue of a complex decomposition
$%
  q=
  (q_0+\widetilde{q}_1f)+
  (\widetilde{q}_2+\widetilde{q}_3f)\phi=
  u_1+u_2\phi
$%
\ ($u_1,\,u_2\in \mathbb{F}$). Note that a corresponding
quaternion $\phi$ for a field $\mathbb{F}$ is not uniquely
determined.

\subsection*{Agreements and notations} In the paper we use
substantially accepted notations. Some notations we define as one
goes along (e.g. $\Lin,\;\Clos$, p.\,\pageref{EqSimplSpectr}).

\textit{"Linear operator"} in $\mathbb{H}$--bimodule means a
\textit{right--side} linear operator. (All the results presented
here hold true if term "bimodule"\ one change by term
"right module". Consideration of precisely bimodules is connected
with the tradition of our previous papers).

\textit{An operator of right multiplication} $R_qx:=xq\;(x\in
H),\;q\in \mathbb{H},$ plays a significant role in studying
operators on quaternion bimodules ($R_q$ is not
$\mathbb{H}$--linear but it is $\mathbb{F}$--linear where the
field $\mathbb{F}$ is generated by $q$; see above). It is useful
to emphasize the following simple properties of an operator of right multiplication
which we use throughout the paper without any
special comments: $R_{p+q}=R_p+R_q;\;R_{pq}=R_qR_p.$\

Note that an $\mathbb{R}$--\textit{linear operator} $A$ is
$\mathbb{F}$--\textit{linear} iff $R_fA=AR_f$ ($f$ is an imaginary
unit generating $\mathbb{F}$). $\mathbb{F}$--\textit{linear
operator} $A$ is $\mathbb{H}$--\textit{linear} iff $R_\phi
A=AR_\phi$.

For further consideration we fix an imaginary unit $f$ and a field
$\mathbb{F}$ generated by $f$. In conclusion note that all
arguments concerning $\mathbb{F}$--modules (orthogonality; etc.)
copy the corresponding ones for $\mathbb{C}$--modules
because of isometrical isomorphism%
\footnote%
{%
   We pay reader's attention to the following remarkable
   phenomenon arising in the quaternion theory of functions and
   operators. This phenomenon is connected with the relation between
   such notions as "equality"\ and "isomorphism". Although every such a field
   $\mathbb{F}$ is isomorphic to $\mathbb{C}$ and therefore it is "not
   interesting"\
   from the \textit{algebraic} point of view nevertheless it is a \textit{unique} "copy"
   \ of $\mathbb{C}$ in $\mathbb{H}$. This fact becomes crucial for \textit{analysis}.
   In particular our theory of differentiability
   of functions of a quaternionic variable (\cite{kt_KarSuchtTyshk}) is based precisely on "playing by"\
   these properties.
}%
\ between $\mathbb{F}$ and $\mathbb{C}$. In this case we use the
corresponding notations $\bot_{\mathbb{F}},\ \oplus_{\mathbb{F}}$;
etc.

\section*{The model of a skew--selfadjoint operator with a simple spectrum}

\subsection*{Spectral Theorem}
The spectral theorem for a skew--selfadjoint operator mentioned in
introduction  can be written as follows.

\begin{theorem}\label{ThSpDecSpPair}
Let $A$ be a skew--selfadjoint operator acting on a Hilbert
quaternion bimodule $H$ (whose domain $\mathcal{D}(A)$ is dense in
$H$), $\mathbb{F}$ be an $\mathbb{R}$--containing subfield of\
$\mathbb{H}$. Then there exists a spectral measure $E$ defined on
the Borel $\sigma$--algebra $\mathfrak{B}(\mathbf{f}_+)$
(consisting of all Borel subsets of the half--axis
$\mathbf{f}_+=\{\tau f\mid \tau\geqslant0\}$), skew--selfadjoint
operator $J$ commuting with $E$ and satisfying the condition
$J^2=-I$  such that the following equalities hold
\begin{align}
   & \mathcal{D}(A)=
         \{
           h\in H\mid\!\!
           \int_{\sigma_{\mathbb{F}_+}(A)}\!\!\!\!\!\!
           |q|^2 \langle E(dq)h,h\rangle<\infty
         \}\,;\label{EqThSpDecSpPair2}\\
   & A=\!\int_{\mathbf{f}_+}\!\!R_{-qf}JE(dq)\,.\label{EqThSpDecSpPair3}
\end{align}
\end{theorem}

The integral in \eqref{EqThSpDecSpPair3} is understood as a strong
limit
$$
  Ah=\lim_{n\rightarrow\infty}
     \!\!\!\!\int_{[0,nf]}\!\!\!\!\!
     R_{-qf}JE(dq)h\,,
     \quad h\in\mathcal{D}(A).
$$

We can determine
\begin{equation}\label{EqMeasE}
   E(\alpha):=
   \begin{cases}
      &\!\!\!\!\!E_\mathbb{F}(\alpha)+E_\mathbb{F}(-\alpha),\, 0\not\in\alpha,\\
      &\!\!\!\!\!E_\mathbb{F}(\alpha)+E_\mathbb{F}(-\alpha)-E_\mathbb{F}(\{0\}),\,0\in\alpha
   \end{cases}
   \quad\big(\alpha\in\mathfrak{B}(\mathbf{f}_+)\big)
\end{equation}
where $E_\mathbb{F}$ is a $\mathbb{F}$--linear spectral measure
defined on Borel subsets of the axis $\mathbf{f}=f\mathbb{R}$; i.e.
the spectral measure of $A$ where $A$ acts on $H$ as on a
$\mathbb{F}$--module. The last is isometrically isomorphic to $H$
as to a $\mathbb{C}$--module, therefore, we can use the classical
spectral theorem for constructing $E_\mathbb{F}$. In addition, $J$
and $E_\mathbb{F}$ are connected by the formula
\begin{equation}\label{EqJDef}
    J=R_f\bigl(E_\mathbb{F}(\mathbf{f}_+)-E_\mathbb{F}(\mathbf{f}_-)\bigr)
\end{equation}
where $\mathbf{f}_-=\{\tau f\mid \tau<0\}$ is the negative
half--axis of $\mathbb{F}$. Concerning \eqref{EqMeasE} note that
$\mathbb{H}$--linearity of $E$ is equivalent to validity of the
equality
\begin{equation}\label{EqEFRphiCommut}
     E_\mathbb{F}(\alpha)R_{\phi}=R_{\phi}E_\mathbb{F}(-\alpha)
     \quad\big(\alpha\in\mathfrak{B}(\mathbf{f}_+)\big).
\end{equation}
In particular,  putting $H^+=E_\mathbb{F}(\mathbf{f}_+)H$ we have
the decomposition
\begin{equation}\label{EqHplusRphiHplus}
  H=H^+\oplus_{\mathbb{F}}R_\phi H^+.
\end{equation}

\subsection*{A skew--selfadjoint operator with a simple spectrum}
Let $A$ be a skew--selfadjoint operator acting on a Hilbert
quaternion bimodule $H$. By theorem~\ref{ThSpDecSpPair} every such
an operator determines the spectral measure $E$ and the
corresponding operator $J$ commuting with $E$.

Denote by $\mathcal{I}$ the set of all (bounded) intervals
$\Delta$ of half-axis
$\mathbf{f}_+;\;\:\Delta=[af,bf];\;a,b\in\mathbb{R};\;0\leqslant
a<b$.

\begin{definition}\label{kt_OpSimplSpectr}
$A$ is called to be an operator \textit{with a simple spectrum} if
there exists a \textit{generating} vector $g\in H$:
\begin{equation}\label{EqSimplSpectr}
   \Clos\Lin\{E(\Delta)g\;|\;\Delta\in\mathcal{I}\}=H
\end{equation}
($\Lin$ is the (right-side) $\mathbb{H}$--span, $\Clos$ is the
$H$--norm closure).
\end{definition}

\subsection*{A special generating vector}
Now we prove the key fact for further consideration  of existence
of a generating vector $g$ with the property
\begin{equation}\label{EqPropSpecGenerVect}
   Jg=R_fg.
\end{equation}

Let $\mathcal{T}=\{T_i\}_{i\in I}$ be a set of operators on $H$,
$g\in H$ is a vector. Define
$$
  \mathcal{C}(\mathcal{T},g)=
  \Clos\Lin\{T_ig\;|\;i\in I\};\quad
  \mathcal{C_{\mathbb{F}}}(\mathcal{T},g)=
  \Clos\Lin_{\mathbb{F}}\{T_ig\;|\;i\in I\}
$$
($\Lin_{\mathbb{F}}$ is the (right--side) $\mathbb{F}$--span). The
following properties are quite elementary:
\begin{align}
   & \mathcal{C}_{\mathbb{F}}(\mathcal{T},g_1+g_2)\subseteq
     \Clos\big(\mathcal{C}_{\mathbb{F}}(\mathcal{T},g_1)+
               \mathcal{C}_{\mathbb{F}}(\mathcal{T},g_2)
          \big);
     \label{EqCsumvectInSumC}\\
   & \mathcal{C}(\mathcal{T},g)=
     \Clos\big(\mathcal{C}_{\mathbb{F}}(\mathcal{T},g)+
               R_\phi\mathcal{C}_{\mathbb{F}}(\mathcal{T},g)
          \big).\label{EqCIsSumCF}
\end{align}
If the set $\mathcal{T}$ is closed respectively to the product of
operators then we have
\begin{equation}\label{EqgInCgIn}
   h\in\mathcal{C}_{\mathbb{F}}(\mathcal{T},g)\Rightarrow
   \mathcal{C}_{\mathbb{F}}(\mathcal{T},h)\subseteq
   \mathcal{C}_{\mathbb{F}}(\mathcal{T},g)\quad
   (\text{\textit{if }}\ \forall\,i,j\in
   I\quad T_iT_j\in\mathcal{T})
\end{equation}
(the same is true for $\mathcal{C}(\cdot,\cdot)$). Denote
$$
  \mathcal{E}=
  \{E(\Delta)\;|\;\Delta\in\mathcal{J}\},\quad
  \mathcal{E}_{\mathbb{F}}=
  \{E_{\mathbb{F}}(\Delta)\;|\;\Delta\in\mathcal{J}\}.
$$
By use of \eqref{EqCsumvectInSumC}~--- \eqref{EqgInCgIn} the following
relations can be easily derived from the properties of a spectral
measure:
\begin{align}
   & \forall\,g\in H\quad
     \mathcal{C}_{\mathbb{F}}(\mathcal{E}_{\mathbb{F}},g)\subseteq
     H^+;\label{EqCEgAddProps1}\\
   & g\in H^+\Rightarrow
     \mathcal{C}(\mathcal{E},g)=
     \mathcal{C}_{\mathbb{F}}(\mathcal{E}_{\mathbb{F}},g);\label{EqCEgAddProps2}\\
   & \forall\,g\in H\quad
     \mathcal{C}(\mathcal{E},R_\phi g)=
     R_\phi\mathcal{C}(\mathcal{E},g);\label{EqCEgAddProps3}\\
   & h\,\bot_{\mathbb{F}}\:\mathcal{C}_{\mathbb{F}}(\mathcal{E}_{\mathbb{F}},g)\Rightarrow
     \mathcal{C}_{\mathbb{F}}(\mathcal{E}_{\mathbb{F}},h+g)=
     \mathcal{C}_{\mathbb{F}}(\mathcal{E}_{\mathbb{F}},h)\oplus_{\mathbb{F}}
     \mathcal{C}_{\mathbb{F}}(\mathcal{E}_{\mathbb{F}},g).\label{EqCEgAddProps4}
\end{align}

\begin{statement}\label{PropSpecGenerVect}
Let $A$ be an operator with a simple spectrum. Then there exists%
\footnote%
{%
  The main idea how to construct a special generating vector in statement~\ref{PropSpecGenerVect}
  is borrowed from the \cite[Prop.\,5.2]{kt_Visv}.
}%
\ a generating vector $g$ satisfying \eqref{EqPropSpecGenerVect}.

\begin{proof}
Let $y$ be an arbitrary generating vector for $A$. Then we have
\begin{equation}\label{EqPropSpecGenerVect1}
   H=\mathcal{C}(\mathcal{E},y).
\end{equation}
By \eqref{EqHplusRphiHplus} for some vectors $y^+,x^+\in H^+\quad
y=y^+ + R_\phi x^+$. Denote
\begin{equation}\label{EqPropSpecGenerVect2}
  \mathfrak{Y}^+=\mathcal{C}_{\mathbb{F}}(\mathcal{E}_{\mathbb{F}},y^+),\quad
  \mathfrak{X}^+=\mathcal{C}_{\mathbb{F}}(\mathcal{E}_{\mathbb{F}},x^+).
\end{equation}
Then we have
\begin{equation}\label{EqPropSpecGenerVect3}
   \mathcal{C}_{\mathbb{F}}(\mathcal{E},y)
   \overset{\eqref{EqCsumvectInSumC}}{\subseteq}
   \Clos\big(
             \mathcal{C}(\mathcal{E},y^+)+
             \mathcal{C}(\mathcal{E},R_\phi x^+)
        \big)
   \overset{\eqref{EqCEgAddProps3},\eqref{EqCEgAddProps2},\eqref{EqPropSpecGenerVect2}}{=}
   \Clos(\mathfrak{Y}^+ + R_\phi\mathfrak{X}^+).
\end{equation}
Hence
\begin{equation}\label{EqPropSpecGenerVect4}
   R_\phi\mathcal{C}_{\mathbb{F}}(\mathcal{E},y)\subseteq
   \Clos(R_\phi\mathfrak{Y}^+ + \mathfrak{X}^+).
\end{equation}
Furthermore,
\begin{equation}\label{EqPropSpecGenerVect5}
   \begin{split}
       &  H\overset{
                    \eqref{EqPropSpecGenerVect1},
                    \eqref{EqCIsSumCF},
                    \eqref{EqPropSpecGenerVect3},
                    \eqref{EqPropSpecGenerVect4}
                   }{\subseteq}
          \Clos\big(
                     \mathfrak{Y}^+ + R_\phi\mathfrak{X}^+ +
                     R_\phi\mathfrak{Y}^+ + \mathfrak{X}^+
               \big)
          \overset{\eqref{EqHplusRphiHplus},\eqref{EqCEgAddProps1}}{=} \\
        & =
          \Clos\big(
                   (\mathfrak{Y}^+ + \mathfrak{X}^+)\oplus_{\mathbb{F}}
                   R_\phi(\mathfrak{Y}^+ + \mathfrak{X}^+)
               \big)=
          \Clos(\mathfrak{Y}^+ + \mathfrak{X}^+)\oplus_{\mathbb{F}}
          R_\phi\Clos(\mathfrak{Y}^+ + \mathfrak{X}^+);
   \end{split}
\end{equation}
then by \eqref{EqHplusRphiHplus}, \eqref{EqCEgAddProps1}
\begin{equation}\label{EqPropSpecGenerVect6}
   H^+=\Clos(\mathfrak{Y}^+ + \mathfrak{X}^+).
\end{equation}
Let $v^+$ be an $\mathbb{F}$--orthogonal projection of the vector
$x^+$ on $H^+\ominus_{\mathbb{F}}\mathfrak{X}^+$; and $x^+=v^+ +
w^+,\;w^+\in\mathfrak{X}^+$. Let $g=y^+ + v^+$. Since $g\in H^+$;
therefore, by \eqref{EqJDef} $g$ satisfies
\eqref{EqPropSpecGenerVect}. Now we can easily prove that $g$ is a
generating vector. Denote
$\mathfrak{V}^+=\mathcal{C}_{\mathbb{F}}(\mathcal{E},v^+)$. By
\eqref{EqgInCgIn}
\begin{equation}\label{EqPropSpecGenerVect7}
   \mathcal{C}_{\mathbb{F}}(\mathcal{E},w^+)\subseteq\mathfrak{Y}^+.
\end{equation}
Hence
\begin{equation}\label{EqPropSpecGenerVect8}
  \mathfrak{Y}^+ + \mathfrak{X}^+=
  \mathfrak{Y}^+ +
  \mathcal{C}_{\mathbb{F}}(\mathcal{E}_{\mathbb{F}},v^+ + w^+)
  \overset{\eqref{EqCEgAddProps4}}{=}
  \mathfrak{Y}^+ +\mathfrak{V}^+ +
  \mathcal{C}_{\mathbb{F}}(\mathcal{E}_{\mathbb{F}},w^+)
  \overset{\eqref{EqPropSpecGenerVect7},\eqref{EqCEgAddProps4}}{=}
  \mathfrak{Y}^+ \oplus_{\mathbb{F}}\mathfrak{V}^+.
\end{equation}
On the other hand
\begin{equation}\label{EqPropSpecGenerVect9}
  H^+
  \overset{\eqref{EqPropSpecGenerVect6},\eqref{EqPropSpecGenerVect8}}{=}
  \mathcal{C}_{\mathbb{F}}(\mathcal{E}_{\mathbb{F}},y^+)\oplus_{\mathbb{F}}
  \mathcal{C}_{\mathbb{F}}(\mathcal{E}_{\mathbb{F}},v^+)
  \overset{\eqref{EqCEgAddProps4}}{=}
  \mathcal{C}_{\mathbb{F}}(\mathcal{E}_{\mathbb{F}},y^+ + v^+)
  \overset{\eqref{EqCEgAddProps2}}{=}
  \mathcal{C}_{\mathbb{F}}(\mathcal{E},g).
\end{equation}
By \eqref{EqPropSpecGenerVect9}, \eqref{EqCIsSumCF},
\eqref{EqHplusRphiHplus} we obtain $H=\mathcal{C}(\mathcal{E},g)$;
i.e. $g$ is a generating vector.
\end{proof}
\end{statement}

\subsection*{Model}
Further we assume that a generating vector $g$ of a
skew--selfadjoint operator $A$ with a simple spectrum satisfies
\eqref{EqPropSpecGenerVect1}.

Consider the operator
$$
  (Qh)(\lambda)=\lambda h(\lambda)
$$
acting on the $\mathbb{H}$--bimodule
$L^2_{\sigma}(\mathbf{f}_+,\mathbb{H})$ of all square integrable
by a measure $\sigma$ $\mathbb{H}$--valued functions on the
half--axis $\mathbf{f}_+$. Since the domain of $Q$ is
$$
  \mathcal{D}(Q)=
  \big\{
        h\in L^2_{\sigma}\;|\;
        \int_{\mathbf{f}_+}\!
        |\lambda|^2\,|h(\lambda)|^2\sigma(d\lambda)<\infty
  \big\};
$$
therefore, by standard arguments one can show that $Q$ is  a
skew--selfadjoint operator.

To find out the spectral measure $E_\mathbb{F}$ we use a slight
modification of the standard algorithm which was considered in our
paper \cite{kt_KarTyshkSpDecUchZ}. In particular, for any interval
$(af,bf),\;-\infty<a<b<\infty$; and any function $h$ from
$\mathcal{D}(Q)$ which has the form
$h(\lambda)=h_1(\lambda)+h_2(\lambda)\phi$ where
$h_1(\lambda),h_2(\lambda)\in \mathbb{F}$ we  have
\begin{align*}
   & E_\mathbb{F}\big((af,bf)\big)h(\lambda)=
     \chi_{(a,b)}(-\lambda f)h_1(\lambda)+\chi_{(a,b)}(\lambda f)h_2(\lambda)\phi
     =\\
   & =
     \chi_{(af,bf)}(\lambda)h_1(\lambda)+\chi_{(-bf,-af)}(\lambda)h_2(\lambda)\phi.
\end{align*}
Then for any $\alpha\in \mathfrak{B}(\mathbf{f})$
\begin{equation*}
   E_\mathbb{F}(\alpha)h(\lambda)=
   \chi_{\alpha}(\lambda)h_1(\lambda)+\chi_{-\alpha}(\lambda)h_2(\lambda)\phi.
\end{equation*}
Hence for $\alpha\in \mathfrak{B}(\mathbf{f}_+\backslash\{0\})$
\begin{equation*}
   E_\mathbb{F}(\alpha)h(\lambda)=
   \chi_{\alpha}(\lambda)h_1(\lambda),\,\,E_\mathbb{F}(-\alpha)h(\lambda)=
   \chi_{\alpha}(\lambda)h_2(\lambda)\phi;
\end{equation*}
and, finally,
$$
  E(\alpha)h(\lambda)= \chi_{\alpha}(\lambda)h(\lambda).
$$

If $\alpha\in \mathfrak{B}(\mathbf{f}_+),\;0\in \alpha$, then
\begin{equation*}
\begin{cases}
   E_\mathbb{F}(\alpha)h(\lambda) & \!\!\!\!=
     \chi_{\alpha}(\lambda)h_1(\lambda)+\widetilde{h}_2(\lambda)\phi,\\
   E_\mathbb{F}(-\alpha)h(\lambda)& \!\!\!\!=
     \widetilde{h}_1(\lambda)+\chi_{\alpha}(\lambda)h_2(\lambda)\phi
\end{cases}
\end{equation*}
where
\begin{equation*}
\widetilde{h}_r(\lambda)=
\begin{cases}
   h_r(0),& \!\!\!\lambda=0,\\
   0,     & \!\!\!\lambda\neq 0,
\end{cases}
\quad(r=1,2).
\end{equation*}
By the equality
\begin{equation*}
E_\mathbb{F}(\{0\})h(\lambda)=
\begin{cases}
   h(0),& \!\!\!\lambda=0,\\
   0,& \!\!\!\lambda\neq 0
\end{cases}
\end{equation*}
one can also obtain the equality
$$
  E(\alpha)h(\lambda)=
  \chi_{\alpha}(\lambda)h(\lambda).
$$
By virtue of \eqref{EqJDef} and above formulae the corresponding
operator $J$ for $Q$ has the form
\begin{equation*}
   (Jh)(\lambda)=
   (h_1(\lambda)-h_2(\lambda)\phi)f.
\end{equation*}
(almost everywhere and except $0$)

Define the step function
\begin{equation*}
   g(\lambda)=\alpha_k,\;\lambda\in
   \Delta_k=
   \big((k-1)f,kf\big),\;\alpha_k\in \mathbf{f}_+\backslash\{0\},
\end{equation*}
with the condition $\sum |\alpha_k|^2\sigma(\Delta_k)<\infty$.
This function generates $Q$ in the sense of
definition~\ref{kt_OpSimplSpectr}. Indeed, the (right) quaternionic
span of the set of functions $E(\Delta)g$ coincides with the set
of all finite step functions and this set is dense in
$L^2_{\sigma}(\mathbf{f}_+,\mathbb{H})$. Note that
$(Jg)(\lambda)=\alpha_k f,\,\,\lambda\in \Delta_k$; i.e.
$Jg=R_fg$. Briefly denote $L^2_{\sigma}(\mathbf{f}_+,\mathbb{H})$
by $L^2_{\sigma}$.

\begin{theorem}\label{ThSkewselfadjOpModelAux}
Let $A$ be a skew--selfadjoint operator with a simple spectrum
acting on a Hilbert quaternion bimodule $H$, $g$ be a generating
vector, $\sigma(\alpha)=\langle E(\alpha)g,g\rangle$ be a measure
defined on the $\sigma$--algebra $\mathfrak{B}(\mathbf{f}_+)$.
Then the map $\Phi:L^2_{\sigma}\rightarrow H$ determined by the
integral
\begin{equation*}
   \Phi h=
   \int_{\mathbf{f}_+}\!\!R_{h(\lambda)}E(d\lambda)g,
\end{equation*}
sets an isometric isomorphism between $L^2_{\sigma}$ and $H$ such
that
\begin{equation*}
   A\Phi h=
   \int_{\mathbf{f}_+}\!\!R_{\lambda h(\lambda)}E(d\lambda)g.
\end{equation*}

\begin{proof}
Denote
\begin{equation*}
   G=
   \big\{
         \widehat{h}\in H\;|\;
         \exists\,h\in L^2_{\sigma}:\;
         \widehat{h}=
         \int_{\mathbf{f}_+}\!\!R_{h(\lambda)}E(d\lambda)g
   \big\}=
   \Im(\Phi)
\end{equation*}
(range of $\Phi$). Let $\Delta$ be an interval of $\mathbf{f}_+$.
Then $\chi_{\Delta}\in L^2_{\sigma}$; and
\begin{equation*}
   \Phi\chi_{\Delta}=
   \int_{\mathbf{f}_+}\!\!E(d\lambda)g\chi_{\Delta}(\lambda)=
   \int_{\Delta}\!\!E(d\lambda)g=
   E(\Delta)g.
\end{equation*}
Hence $G$ is dense in $H$. Let $\widehat{h}\in G.$ Then
$\langle\widehat{h},\widehat{h}\rangle_{H}=\int_{\mathbf{f}_+}\langle
E(d\lambda)g,\widehat{h}\rangle h(\lambda)$. Since
\begin{align*}
   & \langle E(\alpha)g,\widehat{h}\rangle=
     \int_{\mathbf{f}_+}\!\langle E(\alpha)g,E(d\lambda)gh(\lambda)\rangle=
     \int_{\mathbf{f}_+}\!\overline{h(\lambda)}\langle E(d\lambda)E(\alpha)g,g\rangle
     = \\
   & =\int_{\alpha}\!\overline{h(\lambda)}\langle E(d\lambda)g,g\rangle=
     \int_{\alpha}\!\overline{h(\lambda)}\sigma(d\lambda);
\end{align*}
therefore,
\begin{equation*}
   \langle\widehat{h},\widehat{h}\rangle_{H}=
   \int_{\mathbf{f}_+}\!|h(\lambda)|^2\sigma(d\lambda)=
   \langle h(\lambda),h(\lambda)\rangle_{L^2_{\sigma}}.
\end{equation*}
Thus, $\Phi$ isometrically maps the dense subset of $L^2_{\sigma}$
onto a dense subset of $G$. Density of $L^2_{\sigma}$ yields
closedness of $G$. Hence $G=H$; and $\Phi$ is an unitary operator.

Next prove the second part of the theorem. Let $h$ be a finite
function from $L^2_{\sigma}$ with  a support $[af,bf]$;
$\widehat{h}=\int_{\mathbf{f}_+}R_{h(\lambda)}E(d\lambda)g$. To
prove that $\widehat{h}\in D(A)$ we consider the integral
$\int_{\mathbf{f}_+}|\lambda|^2 \langle
E(d\lambda)\widehat{h},\widehat{h}\rangle$. Since
\begin{equation*}
   \langle E(\alpha)\widehat{h},\widehat{h}\rangle=
   \!\int_{\mathbf{f}_+}\!\overline{h(\nu)}\langle E(\alpha)\widehat{h},E(d\nu)g\rangle=
   \!\int_{\mathbf{f}_+}\!\overline{h(\nu)}\langle \widehat{h},E(\alpha)E(d\nu)g\rangle=
   \!\int_{\alpha}\!\overline{h(\nu)}\langle\widehat{h},E(d\nu)g\rangle
\end{equation*}
and
\begin{equation*}
   \langle \widehat{h},E(\beta)g\rangle=
   \!\int_{\mathbf{f}_+}\!\langle E(d\xi)g,E(\beta)g\rangle h(\xi)=
   \!\int_{\beta}\!\langle E(d\xi)g,g\rangle h(\xi)=
   \!\int_{\beta}\!\sigma(d\xi)h(\xi);
\end{equation*}
therefore,
\begin{equation*}
   \langle E(\alpha)\widehat{h},\widehat{h}\rangle=
   \int_{\alpha}\!|h(\nu)|^2\sigma(d\nu).
\end{equation*}
Hence
\begin{equation*}
   \int_{\mathbf{f}_+}\!|\lambda|^2\langle E(d\lambda)\widehat{h},\widehat{h}\rangle=
   \!\int_{\mathbf{f}_+}\!|\lambda|^2|h(\lambda)|^2\sigma(d\lambda)=
   \!\!\int_{[af,bf]}\!\!\!\!|\lambda|^2|h(\lambda)|^2\sigma(d\lambda)<\infty.
\end{equation*}
Further, for any vector $\widehat{h}=\Phi h\in\mathcal{D}(A)\quad
A\widehat{h}=\int_{\mathbf{f}_+}\!\!R_{-\lambda f}J
E(d\lambda)\widehat{h}$. Since
\begin{equation*}
   E(\alpha)\widehat{h}=
   \!\int_{\mathbf{f}_+}\!\!E(\alpha)R_{h(nu)}E(d\nu)g=
   \!\int_{\mathbf{f}_+}\!\!R_{h(\nu)}E(\alpha)E(d\nu)g=
   \!\int_{\alpha}\!\!R_{h(\nu)}E(d\nu)g;
\end{equation*}
therefore, by linearity of $R_{-\lambda f}$ and $J$; and
commutativity of $J$ and $E$ we have
\begin{align*}
   & A\widehat{h}=
     \!\int_{\mathbf{f}_+}\!\!R_{-\lambda f}JR_{h(\lambda)}E(d\lambda)g=
     \!\int_{\mathbf{f}_+}\!\!R_{h(\lambda)}R_{-\lambda f}E(d\lambda)Jg=
     \!\int_{\mathbf{f}_+}\!\!R_{h(\lambda)}R_{-\lambda f}E(d\lambda)R_fg
     =\\
     =
   & \int_{\mathbf{f}_+}\!\!R_{h(\lambda)}R_{\lambda}E(d\lambda)g=
     \!\int_{\mathbf{f}_+}\!\!R_{\lambda h(\lambda)}E(d\lambda)g.
\end{align*}
\end{proof}
\end{theorem}

Thus, on  the set of finite functions from $L^2_{\sigma}$ the equality
$A\Phi=\Phi Q$ holds true (or, in other words, $A=\Phi
Q\Phi^{-1}$). As a consequence we obtain the main result of this
paper.

\begin{theorem}\label{ThSkewselfadjOpModelMain}
Any skew--selfadjoint operator with a simple spectrum acting on a
Hilbert quaternion bimodule is unitarily equivalent to the
operator of left multiplication by an independent variable on a
functional bimodule $L^2_{\sigma}$.
\end{theorem}

\end{document}